\documentclass[12pt]{article}
\usepackage{amsmath}
\usepackage{amssymb}
\usepackage{amscd}
\usepackage{amsthm}

\theoremstyle{plain}
\newtheorem{Thm}{Theorem}[section]
\newtheorem{Conj}[Thm]{Conjecture}
\newtheorem{Prop}[Thm]{Proposition}
\newtheorem{Cor}[Thm]{Corollary}
\newtheorem{Lem}[Thm]{Lemma}

\theoremstyle{definition}
\newtheorem{Defn}[Thm]{Definition}
\newtheorem{Expl}[Thm]{Example}

\newtheorem{Rem}[Thm]{Remark}

\numberwithin{equation}{section}

\title{$D$-equivalence and $K$-equivalence}
\author{Yujiro Kawamata}

\begin{document}

\maketitle

\begin{abstract}
Let $X$ and $Y$ be smooth projective varieties over $\mathbb{C}$.
They are called {\it $D$-equivalent} if their derived categories of 
bounded complexes of coherent sheaves are equivalent as triangulated 
categories, while 
{\it $K$-equivalent} if they are birationally equivalent and the pull-backs of 
their canonical divisors to a common resolution coincide.
We expect that the two equivalences coincide for birationally 
equivalent varieties.
We shall provide a partial answer to the above problem in this paper.
\end{abstract}

\section{Introduction}

Let $X$ be a smooth projective variety.
We denote by $D(X) = D^b(\text{Coh}(X))$ the derived category of 
bounded complexes of coherent sheaves on $X$
(in \S 6, we shall consider a generalization where $X$ has singularities).
It is known that $D(X)$ has a structure of a triangulated category.

\begin{Defn}
Let $X$ and $Y$ be smooth projective varieties.
They are called {\it $D$-equivalent} if their derived categories 
$D(X)$ and $D(Y)$ of 
bounded complexes of coherent sheaves are equivalent
as triangulated categories, i.e., there exists an equivalence
of categories $\Phi: D(X) \to D(Y)$ which commutes with the translations 
and sends any distinguished triangle to a 
distinguished triangle.
They are called {\it $K$-equivalent}
if they are birationally equivalent and if there exists a 
smooth projective variety $Z$ with birational morphisms
$f: Z \to X$ and $g: Z \to Y$ such that the pull-backs of the canonical
divisors are linearly equivalent: $f^*K_X \sim g^*K_Y$.
\end{Defn}

We shall consider the following conjecture which predicts 
that the $D$ and $K$-equivalences
coincide for birationally equivalent varieties.

\begin{Conj}\label{D and K}
Let $X$ and $Y$ be birationally equivalent smooth projective varieties.
Then the following are equivalent.

(1) There exists an equivalence of triangulated categories
$D(X) \cong D(Y)$.

(2) There exists a smooth projective 
variety $Z$ and birational morphisms $f: Z \to X$ and $g: Z \to Y$ such that
$f^*K_X \sim g^*K_Y$.
\end{Conj}

The category of coherent sheaves $\text{Coh}(X)$ reflects the biregular 
geometry of $X$, 
but we expect that the derived category $D(X)$ captures more 
essential properties such as its birational geometry.

A derived category is a purely algebraic object.
But one can sometimes recover the geometry from it:

\begin{Thm}\cite{BO2}\label{BO2}
Let $X$ be a smooth projective variety.
Assume that $K_X$ or $-K_X$ is ample.

(1) Let $Y$ be another smooth projective variety. 
Assume that there exists an equivalence of categories $\Phi: 
D(X) \to D(Y)$ which commutes with the translations.
Then there is an isomorphism $\phi: X \to Y$.

(2) The group of isomorphism classes of exact autoequivalences of $D(X)$ 
is isomorphic to the semi-direct product of $\text{Aut}(X)$ and 
$\text{Pic}(X) \oplus \mathbb{Z}$.
\end{Thm}

We shall prove a generalization of Bondal-Orlov's theorem in this paper:

\begin{Thm}\label{D to K}($=$ Theorem~\ref{D to K'})
Let $X$ and $Y$ be smooth projective varieties.
Assume that the bounded derived categories of coherent sheaves on them are 
equivalent as triangulated categories: $D(X) \cong D(Y)$.
Then the following hold:

(0) $\dim X = \dim Y$. Let $n$ be the common dimension.

(1) If $K_X$ (resp. $-K_X$) is nef, then $K_Y$ (resp. $-K_Y$) is also nef, 
and an equality on the numerical Kodaira dimension
$\nu(X) = \nu(K_Y)$ (resp. $\nu(X, -K_X) = \nu(Y, -K_Y)$) holds.

(2) If $\kappa(X) = n$, i.e., $X$ is of general type, or if 
$\kappa(X, -K_X) = n$, then $X$ and $Y$ are birationally equivalent. 
Moreover, there exist birational morphisms 
$f: Z \to X$ and $g: Z \to Y$ from a 
smooth projective variety $Z$ such that $f^*K_X \sim g^*K_Y$.
\end{Thm}

We also consider the following conjecture:

\begin{Conj}\label{FM}
For a given smooth projective variety $X$, there exist only finitely many 
smooth projective varieties $Y$ up to isomorphisms such that 
$D(Y)$ is equivalent to $D(X)$ as a triangulated category.
\end{Conj}

We shall give an affirmative answer for surfaces in \S 3 by
extending a result of Bridgeland and Maciocia \cite{BM2}:

\begin{Thm}($=$ Theorems~\ref{FM2'} and \ref{FM2})
Let $X$ be a smooth projective surface.
Then there exist at most finitely many smooth projective surfaces $Y$
up to isomorphism such that the derived categories 
$D(X)$ and $D(Y)$ are equivalent as triangulated categories.
Moreover, if $X$ contains a $(-1)$-curve but 
is not isomorphic to a relatively minimal elliptic rational 
surface, then any such $Y$ is isomorphic to $X$.
\end{Thm}

The above conjecture can be regarded as a generalization of 
the conjecture which predicts that there exist only finitely many minimal 
models up to isomorphisms in a fixed birational equivalence class 
(\cite{CY}).
Note that we do not assume the minimality of $X$ in Conjecture~\ref{FM}.

We consider the reverse direction from $K$-equivalence to $D$-equivalence
in the latter half of the paper.
We collects some facts from minimal model theory in \S 4, and 
we calculate some examples in arbitrary dimension in \S 5.
In the case of dimension $3$, we have a complete answer even for 
the case of singular varieties:

\begin{Thm}\label{K to D'}($=$ Theorems~\ref{$3$-flop} and \ref{flop to D})
Let $X$ and $Y$ be normal projective varieties of dimension $3$ 
having only $\mathbb{Q}$-factorial terminal singularities, and 
let $\mathcal{X}$ and $\mathcal{Y}$ be their canonical covering stacks.
Assume that $X$ and $Y$ are $K$-equivalent.
Then the bounded derived categories of coherent orbifold sheaves 
$D(\mathcal{X})$ and $D(\mathcal{Y})$ are equivalent as triangulated 
categories.
\end{Thm}

Acknowledgement: The author would like to thank Tom Bridgeland,
Jiun-Cheng Chen, Akira Ishii, Keiji Oguiso, Burt Totaro and Jan Wierzba 
for useful discussions or comments and the anonymous referee for suggestions.


\section{From $D$-equivalence to $K$-equivalence}

We need the concept of Fourier-Mukai transformation:

\begin{Defn}
Let $X$ and $Y$ be smooth projective varieties, 
and let $p_1: X \times Y \to X$ and $p_2: X \times Y \to Y$ be projections.
For an object $e \in D(X \times Y)$, we define an {\it integral functor}
$\Phi^e_{X \to Y}: D(X) \to D(Y)$ by 
\[
\Phi^e_{X \to Y}(a) = p_{2*}(p_1^*(a) \otimes e)
\]
for $a \in D(X)$, where $p_1^*$ and $\otimes$ are the 
right derived functors and
$p_{2*}$ is the left derived functor.
An integral functor is said to be a {\it Fourier-Mukai transformation} if
it is an equivalence.
\end{Defn}

The following theorem by Orlov is fundamental for the proof of
Theorem~\ref{D to K'}.

\begin{Thm}\cite{O1}\label{O1}
Let $\Phi: D(X) \to D(Y)$ be a functor of bounded derived categories
of coherent sheaves which commutes with the translations 
and sends any distinguished triangle to a 
distinguished triangle.
Assume that $\Phi$ is fully faithful and has a right adjoint.
Then there exists an object $e \in D(X \times Y)$ such that 
$\Phi$ is isomorphic to the integral functor $\Phi^e_{X \to Y}$.
Moreover, $e$ is uniquely determined up to isomorphisms.
\end{Thm}

The following theorem guarantees that the $D$-equivalence implies the
$K$-equivalence at least for general type varieties.

\begin{Thm}\label{D to K'}
Let $X$ and $Y$ be smooth projective varieties.
Assume that the bounded derived categories of coherent sheaves on them are 
equivalent as triangulated categories: $D(X) \cong D(Y)$.
Then the following hold:

(0) $\dim X = \dim Y$. Let $n$ be the common dimension.

(1) If $K_X$ (resp. $-K_X$) is nef, then $K_Y$ (resp. $-K_Y$) is also nef, 
and an equality on the numerical Kodaira dimension
$\nu(X) = \nu(Y)$ (resp. $\nu(X, -K_X) = \nu(Y, -K_Y)$) holds.

(2) If $\kappa(X) = n$, i.e., $X$ is of general type, or if 
$\kappa(X, -K_X) = n$, then $X$ and $Y$ are birationally equivalent. 
Moreover, there exist birational morphisms 
$f: Z \to X$ and $g: Z \to Y$ from a 
smooth projective variety $Z$ such that $f^*K_X \sim g^*K_Y$.
\end{Thm}

\begin{proof}
By Theorem~\ref{O1}, there exists an object $e \in D(X \times Y)$ 
such that $\Phi^e_{X \to Y}: D(X) \to D(Y)$ is an equivalence.
Let 
\[
e^{\vee} = RHom_{\mathcal{O}_{X \times Y}}(e, \mathcal{O}_{X \times Y})
\]
the the derived dual object.
By the Grothendieck duality, 
the right and left adjoint functors of $\Phi=\Phi^e_{X \to Y}$ are given by
$\Phi^{e^{\vee} \otimes p_1^*\omega_X[\dim X]}_{Y \to X}$ and 
$\Phi^{e^{\vee} \otimes p_2^*\omega_Y[\dim Y]}_{Y \to X}$.

Since $\Phi$ is an equivalence, the right and left adjoint functors of 
$\Phi=\Phi^e_{X \to Y}$ are isomorphic.
By Theorem~\ref{O1} again, we have an isomorphism of objects
\[
e^{\vee} \otimes p_1^*\omega_X[\dim X]
\cong e^{\vee} \otimes p_2^*\omega_Y[\dim Y].
\]
It follows immediately that $\dim X = \dim Y$.

Let $H^i(e^{\vee})$ be the cohomology sheaves, 
$\Gamma$ the union of the 
supports of the $H^i(e^{\vee})$ for all $i$, 
$\Gamma = \bigcup_j Z_j$ the decomposition to irreducible components,
and let $\nu_j : \tilde Z_j \to Z_j$ be the normalizations.
We take a $Z_j$ and assume that it is an irreducible component of the
support of $H^i(e^{\vee})$.
By taking the determinant of both sides of the isomorphism 
\[
\nu_j^*(H^i(e^{\vee}) \otimes p_1^*\omega_X) \cong 
\nu_j^*(H^i(e^{\vee}) \otimes p_2^*\omega_Y)
\] 
we obtain
\[
\nu_j^*p_1^*\omega_X^{\otimes m_j} \cong \nu_j^*p_2^*\omega_Y^{\otimes m_j}
\]
where $m_j$ is the rank of $\nu_j^*H^i(e^{\vee})$.

(1) Since $\Phi^e_{X \to Y}$ is an equivalence, 
the projections $p_1 \vert_{\Gamma}: \Gamma \to X$
and $p_2 \vert_{\Gamma}: \Gamma \to Y$ are surjective.
Let $Z_1$ be an irreducible component of $\Gamma$ which dominates $Y$.
If $K_X$ is nef, then 
$m_1\nu_1^*p_1^*K_X \sim m_1\nu_1^*p_2^*K_Y$ 
is also nef, hence so is $K_Y$.  
We have also 
$\nu(X) \ge \nu(\tilde Z_1, \nu_1^*p_2^*K_Y) = \nu(Y)$, thus 
$\nu(X) = \nu(Y)$.
The case where $-K_X$ is nef is proved similarly.

(2) If $\kappa(X) = n$, then there exist an ample $\mathbb{Q}$-divisor $A$ 
and an effective $\mathbb{Q}$-divisor $B$ on $X$ 
such that $K_X \sim_{\mathbb{Q}} A + B$ by Kodaira's lemma.
Let $Z_1$ be an irreducible component of $\Gamma$ which dominates $X$.
Then the projection 
$p_2 \vert_{Z_1}: Z_1 \to Y$ 
is quasi-finite on $Z_1 \setminus p_1^{-1}(\text{Supp}(B))$.
Indeed, if there exists a curve $C$ which is contained in 
$Z_1 \cap p_2^{-1}(y)$ for a point $y \in Y$ but not entirely in 
$p_1^{-1}(\text{Supp}(B))$, then we have
$(p_2^*K_Y \cdot C) = 0$ while
$(p_1^*K_X \cdot C) \ge (p_1^*A \cdot C) > 0$, a contradiction.
Since $\dim X = \dim Y = n$, it follows that $\dim Z_1 = n$ and 
$Z_1$ also dominates $Y$.

We claim that the set $\Gamma \cap p_1^{-1}(x)$ consisits
of $1$ point for a general point $x \in X$.
Indeed, the previous argument showed already that $\Gamma \cap p_1^{-1}(x)$ is 
a finite set.
If it is not connected, then the natural map
$\text{Hom}_{D(X)}(\mathcal{O}_x, \mathcal{O}_x) \to 
\text{Hom}_{D(Y)}(\Phi(\mathcal{O}_x), \Phi(\mathcal{O}_x))$
is not surjective, a contradiction.
Therefore, $Z_1$ is a graph of a birational map.
If we take $Z$ to be any resolution of $Z_1$, then the conclusion holds.

The case where $\kappa(X, -K_X)=n$ is proved similarly. 
\end{proof}

\begin{Rem}
(0) The differential geometric picture of the above proof is that 
the kernel object $e$ of the Fourier-Mukai transformation 
cannot spread itself if the Ricci curvature is non-vanishing.

(1) In the case where $K_X$ or $-K_X$ is ample, 
we can also reprove Theorem~\ref{BO2}~(2) 
by a similar argument as above.

Indeed, if we take $B=0$, then $Z_1$ becomes a graph of an isomorphism,
say $h$.
Now $e$ can be considered as a complex of sheaves on $X$ so that
we have
$\Phi(\mathcal{O}_x) \cong h(e \otimes_{\mathcal{O}_X} \mathcal{O}_x)$
for any $x \in X$,
where the tensor product is taken in $D(X)$.
Since 
\[
\text{Hom}_{D(X)}^p(\Phi(\mathcal{O}_x), \Phi(\mathcal{O}_x)) = 0
\]
for any $p < 0$, it follows that there exists an integer $i_0$ 
such that $e[i_0]$ is a sheaf.
Since 
\[
\text{Hom}_{D(X)}(\Phi(\mathcal{O}_x), \Phi(\mathcal{O}_x)) = \mathbb{C}
\]
$e[i_0]$ is invertible.

We note that we did not assume in Theorem~\ref{BO2} 
that $\Phi$ sends any distinguished triangle to a 
distinguished triangle.

(2) We can extend Theorem~\ref{BO2}~(2) to the case where $X$ admits 
quotient singularities if $K_X$ generates the local class group
at any point as in \cite{Francia}.
Namely, let $\mathcal{X}$ be the smooth stack which lies naturally above $X$ 
and let $D(\mathcal{X}) = D^b(\text{Coh}(\mathcal{X}))$ 
be the derived category of bounded complexes of coherent sheaves 
on $\mathcal{X}$ (see \S 6).
Then $\text{Auteq}(D(\mathcal{X}))$
is isomorphic to the semi-direct product of $\text{Aut}(X)$ and 
$\text{Pic}(\mathcal{X}) \oplus \mathbb{Z}$.
The proof is the same as in \cite{BO2}.

On the other hand, if $K_X$ does not generate the local class group, then
the group of autoequivalences is much larger.
For example, if $Y$ is a smooth projective minimal surface of general 
type and $X$ is its canonical model, then 
$D(\mathcal{X})$ is equivalent to $D(Y)$.
If $C$ is an exceptional curve of the resolution $Y \to X$, then
$\mathcal{O}_C(-1)$ is a $2$-spherical object in $D(Y)$ and generates an 
autoequivalence of infinite order (\cite{ST}, see also \S 4).

(3) If $\nu(X) = \nu(Y) = 0$ in Theorem~\ref{D to K'}~(1), then 
$K_X \sim 0$ if and only if $K_Y \sim 0$ because $\Phi$ commutes with
the Serre functors. 
More generally, it is known that 
the orders of the canonical divisors coincide (\cite{BM2}~Lemma~2.1).
\end{Rem}


\section{Fourier-Mukai partners of surfaces}

We have a complete picture of $D$ and $K$-equivalences for
surfaces. 
We start with the case of minimal surfaces:

\begin{Thm}\cite{BM2}\label{FM2'}
Let $X$ be a smooth projective surface.
Assume that there is no $(-1)$-curve on $X$.
Then there exist at most finitely many smooth projective surfaces $Y$
such that the derived categories 
$D(X)$ and $D(Y)$ are equivalent as triangulated categories.
\end{Thm}

We note that there are Fourier-Mukai partners which are not birationally
equivalent in the 
case of abelian or K3 or elliptic surfaces (\cite{M1}, \cite{M2},
\cite{O1}, \cite{O2}, \cite{BM2}, \cite{HLOY}).
It is rather surprising that the existence of a $(-1)$-curve reduces
the symmetry drastically:

\begin{Thm}\label{FM2}
Let $X$ be a smooth projective surface.
Assume that there exists a $(-1)$-curve on $X$.
Then there exist at most finitely many smooth projective surfaces $Y$
such that the derived categories 
$D(X)$ and $D(Y)$ are equivalent as triangulated categories.
Moreover, if $X$ is not isomorphic to a relatively minimal elliptic rational 
surface, then any such $Y$ is isomorphic to $X$.
\end{Thm}

\begin{proof}
We use the notation of the proof of Theorem~\ref{D to K'}.
Let $C$ be a $(-1)$-curve and $\Gamma_C = p_1^{-1}(C) \cap \Gamma$. 
Since $- K_X \vert_C$ is ample, the projection 
$p_2 \vert_{\Gamma_C}: \Gamma_C \to Y$ is a finite morphism. 
We have two possibilities that $\dim \Gamma_C = 1$ or $2$.

Assume first that $\dim \Gamma_C = 1$.
We take an irreducible component $Z_1$ of $\Gamma$ which dominates $X$,
and let $Z_{1,C} = p_1^{-1}(C) \cap Z_1$ and $C' = p_2(Z_{1,C})$.
We know that $\dim Z_{1,C} = \dim C' = 1$.
It follows that $\dim Z_1 = 2$ and the projection 
$p_1 \vert_{Z_1}: Z_1 \to X$ is generically finite, hence
a birational morphism as in the proof of Theorem~\ref{D to K'}.

If $Z_1$ dominates $Y$, then the other projection 
$p_2 \vert_{Z_1}: Z_1 \to Y$ is also birational, 
and $X$ and $Y$ are $K$-equivalent through $Z_1$.
Hence $X$ and $Y$ are isomorphic (cf. Lemma~\ref{EF}).

Otherwise, we have $p_2(Z_1) = C'$.
There exists an open dense subset $U \subset X$ such that 
$p_1$ induces an isomorphism 
$p_1^{-1}(U) \cap \Gamma = p_1^{-1}(U) \cap Z_1 \to U$.
Take two distinct points $x_1, x_2 \in U$ which correspond to the 
same point $y \in C'$, i.e., 
$y = p_2(p_1^{-1}(x_1) \cap \Gamma) = p_2(p_1^{-1}(x_2) \cap \Gamma)$.
Then both $\Phi(\mathcal{O}_{x_1})$ and $\Phi(\mathcal{O}_{x_2})$ are
supported at $y$, hence
$\text{Hom}_{D(Y)}^p(\Phi(\mathcal{O}_{x_1}),\Phi(\mathcal{O}_{x_2})) \ne 0$
for some $p$, a contradiction.

Assume next that $\dim \Gamma_C = 2$.
Then $p_2 \vert_{\Gamma_C}: \Gamma_C \to Y$ is dominant.
Since $(K_X \cdot C) < 0$, we deduce that 
$- K_Y$ is nef and $\nu(Y, -K_Y) = 1$.
Hence $- K_X$ is also nef and $\nu(X, -K_X) = 1$ by Theorem~\ref{D to K'}.
By the classification of surfaces, such a surface is isomorphic to 
either a minimal elliptic ruled surface or 
a rational surface with Euler number $12$.
Since $X$ has a $(-1)$-curve, $X$ is a rational surface.
By \cite{BM2}~Proposition~2.3, $Y$ is also a rational surface.

We have the possibilities that $\dim \Gamma = 2$ or $3$.
If $\dim \Gamma = 2$, then we obtain our result as before.
If $\dim \Gamma = 3$, then $X$ and $Y$ are dominated by families of 
curves whose intersection
numbers with the canonical divisors vanish.
Thus $X$ and $Y$ are relatively minimal rational elliptic surfaces.
By \cite{BM2}~Proposition~4.4, we obtain our result.
Here we note that the proof there works also for relatively minimal 
elliptic surfaces of negative Kodaira dimension.
\end{proof}

We can extend some of the above argument to higher dimensional case:

\begin{Prop}\label{blowup}
Let $X$ and $Y$ be smooth projective varieties.
Assume that $\kappa(X) \ge 0$ but $K_X$ is not nef,
and that there is an extremal contraction morphism 
$\phi: X \to W$ which contracts a prime divisor $D$ to a point.
Assume that the derived categories 
$D(X)$ and $D(Y)$ are equivalent as triangulated categories.
Then $X$ and $Y$ are birational and $K$-equivalent.
\end{Prop}

\begin{proof}
We use the notation of the proof of Theorem~\ref{D to K'}.
The proof is similar to that of Theorem~\ref{FM2}.

Let $\Gamma_D = p_1^{-1}(D) \cap \Gamma$. 
If $\dim \Gamma_D = n-1$ for $n = \dim X$, 
then there exists an irreducible component $Z_1$ of $\Gamma$ of 
dimension $n$ which dominates $X$.
Then it follows that $X$ and $Y$ are birational and $K$-equivalent as in
the proof of Theorem~\ref{FM2}.

Assume that $\dim \Gamma_D \ge n$.
Since $- K_X \vert_D$ is ample, the projection 
$p_2 \vert_{\Gamma_D}: \Gamma_D \to Y$ is a finite morphism.
Hence $\dim \Gamma_D = n$, and $- K_Y$ is nef with $\nu(Y, - K_Y) = n-1$,
a contradiction to $\kappa(X)$.
\end{proof}

\begin{Rem}
We cannot expect similar statements for other types of contractions.
For example, let $A$ be an abelian surface, $\hat{A}$ its dual, and $S$ 
a smooth projective surface which contains a $(-1)$-curve.
Let $X = A \times S$ and $Y = \hat{A} \times S$.
Then $X$ has a divisorial contraction, $D(X) \cong D(Y)$, but $X$ and $Y$ are
not birational in general.
\end{Rem}


\section{Flops and minimal models}

We consider normal varieties which are not necessarily smooth in this section.

\begin{Defn}
Let $X$ and $Y$ be normal quasiprojective varieties 
whose canonical divisors are $\mathbb{Q}$-Cartier divisors.
A birational map $\alpha: X -\to Y$ is said to be {\it crepant}
if there exists a smooth quasiprojective variety $Z$ with birational 
projective morphisms
$f: Z \to X$ and $g: Z \to Y$ such that $\alpha \circ f = g$ and 
$f^*K_X \sim_{\mathbb{Q}} g^*K_Y$.
\end{Defn}

\begin{Lem}\label{EF}
Let $\alpha: X -\to Y$ be a crepant birational map between 
quasiprojective varieties with only terminal singularities.
Then $\alpha$ is an isomorphism in codimension $1$; i.e.,  
there exist closed subvarieties $E \subset X$ and $F \subset Y$ 
of codimension at least $2$ such that $\alpha$ induces an isomorphism
$X \setminus E \cong Y \setminus F$.
\end{Lem}

\begin{proof}
Since $X$ has only terminal singularities, a prime divisor $D$ on $Z$ is 
mapped by $f$ to a subvariety of codimension at least $2$ on $X$ if and only if
it appears in the relative canonical divisor $K_{Z/X} = K_Z - f^*K_X$ as an 
irreducible component.
Since a similar statement holds for $g$, our assertion follows from 
the equality $K_{Z/X} = K_{Z/Y}$.
\end{proof}

\begin{Defn}
A projective variety $X$ with only canonical singularities is called
{\it minimal} if $K_X$ is nef.
\end{Defn}

The minimality of a variety is characterized by the minimality of its 
canonical divisor:

\begin{Lem}
Let $X$ and $Y$ be normal projective varieties whose canonical divisors are 
$\mathbb{Q}$-Cartier divisors.
Assume that $X$ and $Y$ are birationally equivalent,
$X$ has only canonical singularities and that $K_X$ is nef.
Then the inequality $K_X \le K_Y$ holds in the following sense:
Let $Z$ any smooth projective variety with projective birational morphisms 
$f: Z \to X$ and $g: Z \to Y$. 
Then there exists a positive integer $m$ such that 
$m(g^*K_Y - f^*K_X)$ is linearly equivalent to an effective divisor.
In particular, any birational map between minimal varieties is crepant. 
\end{Lem}

\begin{proof}
We write $f^*K_X + A = g^*K_Y + B$, where
$A$ and $B$ are effective divisors without common irreducible components.
Since $X$ has only canonical singularities,
we may assume that $\text{codim }g(\text{Supp}(B)) \ge 2$.

Assuming that $B \ne 0$, we shall derive a contradiction.
Let $H$ and $M$ be very ample divisors on $Y$ and $Z$, respectively,
and let $n = \dim Y$ and $d = \dim g(\text{Supp}(B))$.
We consider a generic surface section 
\[
S = g^*H_1 \cap \dots \cap g^*H_d \cap 
M_1 \cap \dots \cap M_{n-d-2}
\]
for $H_i \in \vert H \vert$ and $M_j \in \vert M \vert$.
By the Hodge index theorem, we have $(g^*H^d \cdot M^{n-d-2} \cdot B^2) < 0$,
while $(g^*H^d \cdot M^{n-d-2} \cdot B \cdot (f^*K_X + A - g^*K_Y)) \ge 0$
because $K_X$ is nef and $(g^*H^d \cdot g^*K_Y \cdot M^{n-d-2} \cdot B) = 0$, 
a contradiction.
\end{proof}

We consider a special kind of crepant birational maps called flops:

\begin{Defn}
Let $X$ and $Y$ be quasiprojective varieties with only 
canonical singularities, and $D$ a $\mathbb{Q}$-Cartier divisor on $X$.
A birational map $\alpha: X -\to Y$ is said to be a {\it $D$-flop}, 
or simply a {\it flop}, 
if there exist a normal quasiprojective variety $W$ and crepant birational 
projective morphisms $\phi: X \to W$ and $\psi: Y \to W$ 
which satisfy the following conditions:

(1) $\phi = \psi \circ \alpha$. 

(2) $\phi$ and $\psi$ are isomorphisms in codimension $1$.

(3) $D$ is $\phi$-ample, and for any $\mathbb{Q}$-Cartier divisor $A$ on $X$, 
there exist a $\mathbb{Q}$-Cartier divisor $A_0$ on $W$ and a rational number 
$r$ such that 
$A \sim_{\mathbb{Q}} \phi^*A_0 + rD$.

(4) Let $D'$ be the strict transform of $D$ on $Y$.  Then 
$-D'$ is $\psi$-ample, and for any $\mathbb{Q}$-Cartier divisor $B$ on $Y$, 
there exist a $\mathbb{Q}$-Cartier divisor $B_0$ on $W$ and a rational number 
$r'$ such that 
$B \sim_{\mathbb{Q}} \psi^*B_0 + r'D'$.
\end{Defn}

We can define flops of complex analytic spaces instead of 
quasiprojective varieties in a similar way.
In this case, $X$ and $Y$ are complex analytic spaces which are 
relatively projective over a complex analytic space $W$.

Any crepant birational map between projective varieties with only 
$\mathbb{Q}$-factorial terminal singularities is expected to be 
decomposed into a sequence of flops:

\begin{Thm}\label{$3$-flop}
Let $\alpha: X -\to Y$ be a crepant birational map between 
projective varieties of dimension $3$ 
with only $\mathbb{Q}$-factorial terminal singularities.
Then $\alpha$ is decomposed into a sequence of flops.
\end{Thm}

\begin{proof}
We may assume that the subvariety $E$ of Lemma~\ref{EF} is 
purely $1$-dimensional.
We may also assume that any irreducible component of $E$ is the image of a 
curve on $Z$ which is mapped to a point on $Y$.
Since $\alpha$ is crepant, we have $K_X \vert_E \sim_{\mathbb{Q}} 0$. 
Let $H$ be an ample Cartier divisor on $Y$ such that 
$H - K_Y$ is still ample, and let $H'$ be its strict transform on $X$.
By construction, any curve $C$ such that $(H' \cdot C) \le 0$ is contained in
$E$.
We run the minimal model program with respect to $K_X + \epsilon H'$, where 
$\epsilon$ is a small positive number, for only those extremal rays on which 
$H'$ is non-positive.
Then the associated 
extremal curves are contained in $E$, so we obtain an $H'$-flop.
We denote the result after the flop again by the same letters such as 
$X, E$ and $H'$.
After a finite flops, we have no more extremal rays on which 
$H'$ is non-positive.
Then $H'$ becomes nef and big.
Since $H'$ is ample outside $E$, $H' - K_X$ is also nef and big, 
By the base point free theorem, we obtain a birational morphism
$X \to Y$, which should be an isomorphism. 
\end{proof}


\section{From $K$-equivalence to $D$-equivalence}

The following is a special case of the implication from (2) to (1) in 
Conjecture~\ref{D and K}:

\begin{Conj}\label{F to D}
Let $X$ and $Y$ be smooth projective varieties and 
$\alpha: X \to W \leftarrow Y$ a flop. 
Then there exists an equivalence of triangulated categories
$\Phi: D(X) \to D(Y)$.
\end{Conj}

The examples in this section suggest that
the integral functor $\Phi^e_{X \to Y}$ for the structure sheaf 
$e = \mathcal{O}_{X \times_W Y}$ of the subscheme
$X \times_W Y \subset X \times Y$ might work.

We consider the following $2$ examples of flops in this section.

\begin{Expl}\label{flop}
(1) A {\it standard flop}.
Let $X$ be a smooth projective variety of dimension $2m+1$ 
for some positive integer $m$, and $E$ a subvariety of $X$.  
Assume that $E \cong \mathbb{P}^m$, 
and $N_{E/X} \cong \mathcal{O}_{\mathbb{P}^m}(-1)^{m+1}$.
Let $f: Z \to X$ be the blowing-up with center $E$.
Then the exceptional divisor $G$ is isomorphic to 
$\mathbb{P}^m \times \mathbb{P}^m$ and can be blown-down to another direction,
so that we obtain a birational morphism $g: Z \to Y$ and 
a subvariety $F = g(G) \cong \mathbb{P}^m$.
There is a projective variety $W$ with 
contraction morphisms $\phi: X \to W$ and $\psi: Y \to W$ whose 
exceptional loci are $E$ and $F$, respectively, 
and such that $w_0 = \phi(E) = \psi(F)$ is the only singular point of $W$.
Then $\alpha = g \circ f^{-1} = \psi^{-1} \circ \phi$ is a flop.

(2) {\it Mukai's flop}.
Let $W_0$ be a generic hypersurface section of $W$ in (1) 
through the singular point $w_0$.
Let $X_0 = \phi^{-1}(W_0)$, $Y_0 = \psi^{-1}(W_0)$, 
$\phi_0 = \phi \vert_{X_0}$, and $\psi_0 = \psi \vert_{Y_0}$.
Then $X_0$ and $Y_0$ are smooth, 
and $\alpha_0 = \psi_0^{-1} \circ \phi_0$ is a flop. 
The inverse image $\tilde Z_0 = f^{-1}(X_0) = g^{-1}(Y_0)$ is reducible
with $2$ irreducible components $G$ and $Z_0$, where $Z_0$ is smooth.
The restrictions $f_0 = f \vert_{Z_0}$ and $g_0 = g \vert_{Z_0}$ are 
again birational morphisms, and $\alpha_0 = g_0 \circ f_0^{-1}$.
We set $G_0 = G \cap Z_0$.
Then $f_0(G_0) = E$ and $g_0(G_0) = F$. 
\end{Expl}

We need the following concepts:

\begin{Defn}
A set $\Omega$ of objects of $D(X)$ is said to a {\it spanning class}
if the following hold for any $a \in D(X)$.

(1) $\text{Hom}^p(a, \omega) = 0$ for all $p \in \mathbb{Z}$ and 
all $\omega \in \Omega$ implies that $a \cong 0$

(2) $\text{Hom}^p(\omega,a) = 0$ for all $p \in \mathbb{Z}$ and 
all $\omega \in \Omega$ implies that $a \cong 0$.
\end{Defn}

For example, the set of point sheaves $\{\mathcal{O}_P\}$ 
for a smooth projective variety is a spanning class (\cite{B1}~Example~2.2).

\begin{Defn}
A {\it Serre functor} $S_X: D(X) \to D(X)$ is an autoequivalence of 
triangulated categories which induces bifunctorial isomorphisms
\[
\text{Hom}_{D(X)}(a,b) \to \text{Hom}_{D(X)}(b, S_X(a))^*
\]
for $a,b \in D(X)$.
\end{Defn}

If a Serre functor exists, then it is unique up to isomorphisms.
If $X$ is smooth and projective, then 
$S_X(a) = a \otimes \omega_X[\dim X]$ is a Serre functor.

In order to prove that a functor $\Phi: D(X) \to D(Y)$ to be fully faithful, 
it is sufficient to check it for the spanning class (\cite{B1}~Theorem~2.3):
\[
\Phi: \text{Hom}^p(\omega_1, \omega_2) \cong 
\text{Hom}^p(\Phi(\omega_1), \Phi(\omega_2))
\]
for all $p \in \mathbb{Z}$ and all $\omega_1, \omega_2 \in \Omega$. 
Moreover, by \cite{BKR}~Theorem 2.3, provided that $\Phi = \Phi^e_{X \to Y}$ 
is fully faithful, it is an equivalence if and only if it commutes with the 
Serre functor.
Theorefore, in order to prove our conjecture, 
we may consider locally over an analytic neighborhood of a point of $W$ and 
replace the given flop by any other flop which is 
analytically isomorphic to the original one. 
If $\Phi$ is proved to be fully faithful, 
then it is automatically an equivalence in our case.

\begin{Prop}\cite{BO1}
In Example~\ref{flop}~(1), 
$Z$ is isomorphic to the fiber product $X \times_W Y$ which 
is a closed subscheme of $X \times Y$,
and the functor
\[
g_*f^* = \Phi^{\mathcal{O}_Z}_{X \to Y}: D(X) \to D(Y)
\]
is an equivalence of triangulated categories.
\end{Prop}

\begin{proof}
We may replace $X$, $Y$ and $Z$ by the total space of the 
vector bundles $N_{E/X}$, $N_{F/Y}$ and $N_{G/Z}$, respectively.
We denote by $\mathcal{O}_X(k)$, $\mathcal{O}_Y(l)$ and 
$\mathcal{O}_Z(k,l)$ the pull-backs of $\mathcal{O}_E(k)$,
$\mathcal{O}_F(l)$ and $\mathcal{O}_G(k,l)$, respectively.
The set of objects 
\[
\{\mathcal{O}_X(-k) \in D(X) \vert k=0,1,\dots,m\}
\]
spans $D(X)$.
Since $K_{Z/X} \sim mG$, we have
\[
\begin{CD}
\mathcal{O}_X(-k) @>{f^*}>> \mathcal{O}_Z(-k,0) \cong 
\mathcal{O}_Z(0,k)(kG) @>{g_*}>> \mathcal{O}_Y(k).
\end{CD}
\]
We have 
\[
\text{Hom}^p(\mathcal{O}_X(-k_1), \mathcal{O}_X(-k_1)) \cong 
\text{Hom}^p(\mathcal{O}_Y(k_1), \mathcal{O}_Y(k_1)) \cong 0
\]
for $p \ne 0$ and $k_1, k_2 =0,1,\dots,m$ by the vanishing theorem, and 
\[
\Phi^{\mathcal{O}_Z}_{X \to Y}: 
\text{Hom}(\mathcal{O}_X(-k_1), \mathcal{O}_X(-k_1)) \cong 
\text{Hom}(\mathcal{O}_Y(k_1), \mathcal{O}_Y(k_1))
\]
because $X$ and $Y$ are isomorphic in codimension $1$.
Therefore, $\Phi^{\mathcal{O}_Z}_{X \to Y}$ is an equivalence by the remarks 
preceding to the proposition.
\end{proof}

\begin{Lem}\cite{C}
Let $\pi_X: X \to S$ and $\pi_Y: Y \to S$ be smooth projective morphisms
from smooth quasiprojective varieties
to a smooth quasiprojective curve.
Let $s_0 \in S$ be a point, and let $X_0 = \pi_X^{-1}(s_0)$ and
$Y_0 = \pi_X^{-1}(s_0)$ be fibers.
Let $i_{X_0}: X_0 \to X$, $i_{Y_0}: Y_0 \to Y$ and 
$i_{X \times_S Y}: X \times_S Y \to X \times Y$ be the embeddings.
Let $e \in D(X \times_S Y)$ be an object, and let 
$e_0 = e \otimes \mathcal{O}_{X_0 \times Y_0}$
and $e' = i_{X \times_S Y*}(e)$.
Then there is an isomorphism of functors from $D(X_0)$ to $D(Y)$:
\[
i_{Y_0*} \circ \Phi_{X_0 \to Y_0}^{e_0} 
\cong \Phi_{X \to Y}^{e'} \circ i_{X_0*}.
\]
\end{Lem}

\begin{proof}
Let $i_{X_0 \times Y_0}: X_0 \times Y_0 \to X \times_S Y$ be the embedding.
Let $p_1: X \times_S Y \to X$, $p_2: X \times_S Y \to Y$,
$p_{1,0}: X_0 \times Y_0 \to X_0$ and 
$p_{2,0}: X_0 \times Y_0 \to Y_0$ be projections.
For $a \in D(X_0)$, we have
\[
\begin{split}
&i_{Y_0*} \circ \Phi_{X_0 \to Y_0}^{e_0}(a)
\cong i_{Y_0*}p_{2,0*}(p_{1,0}^*(a) \otimes e_0)
\cong p_{2*}i_{X_0 \times Y_0*}(p_{1,0}^*(a) \otimes e_0) \\
&\cong p_{2*}(i_{X_0 \times Y_0*}p_{1,0}^*(a) \otimes e)
\cong p_{2*}(p_1^*i_{X_0*}(a) \otimes e)
\cong \Phi_{X \to Y}^{e'} \circ i_{X_0*}(a).
\end{split}
\]
\end{proof}

\begin{Cor}\label{Mu1}
In Example~\ref{flop}~(2), the functor
\[
\Phi^{\mathcal{O}_{\tilde Z_0}}_{X_0 \to Y_0}: D(X_0) \to D(Y_0)
\]
is an equivalence of triangulated categories.
\end{Cor}

\begin{proof}
Since $Z = X \times_W Y$ is a subscheme of $X \times_S Y$, we have
the following isomorphisms
\[
i_{X_0*}\Phi^{\mathcal{O}_{\tilde Z_0}(mG)}_{Y_0 \to X_0}
\Phi^{\mathcal{O}_{\tilde Z_0}}_{X_0 \to Y_0}
\cong \Phi^{\mathcal{O}_Z(mG)}_{Y \to X}
\Phi^{\mathcal{O}_Z}_{X \to Y}i_{X_0*}
\cong i_{X_0*}.
\]
For any $a \in D(X_0)$, let $b \in D(X_0)$ be the cone of the natural morphism
\[
a \to \Phi^{\mathcal{O}_{\tilde Z_0}(mG)}_{Y_0 \to X_0}
\Phi^{\mathcal{O}_{\tilde Z_0}}_{X_0 \to Y_0}(a).
\]
Then $i_{X_0*}(b) \cong 0$, hence $b \cong 0$.
\end{proof}

The following concept 
is useful for constructing autoequivalences of derived categories.

\begin{Defn}\cite{ST}
An object $s \in D(X)$ is called {\it $n$-spherical} if 
\[
\text{Hom}^p_{D(X)}(s,s) \cong \begin{cases} \mathbb{C} &\text{ if } p = 0,n \\
0 &\text{ otherwise}. 
\end{cases}
\]
The {\it twisting functors} $T_s, T'_s : D(X) \to D(X)$ are defined 
such that the following triangles are distinguished:
\[
\begin{split}
&R\text{Hom}_X(s,a) \otimes s \to a \to T_s(a) \to 
R\text{Hom}_X(s,a) \otimes s[1] \\
&T'_s(a) \to a \to R\text{Hom}_{\mathbb{C}}(R\text{Hom}_X(a,s), s) \to
T'_s(a)[1]
\end{split}
\]
where $R\text{Hom}_X$ denotes the derived global Hom.
If $s$ is $n$-spherical for $n = \dim X$, then $T_s$ and $T'_s$ are 
equivalences and $T_s \circ T'_s \cong \text{Id}_{D(X)}$.
\end{Defn}

\begin{Expl}
(1) $\mathcal{O}_E$ in Example~\ref{flop}~(1) is a $(2m+1)$-sherical object.
Indeed, since $N_{E/X} \cong \mathcal{O}_E(-1)^{m+1}$, we have
\[
\mathcal{E}xt^p_{\mathcal{O}_{X_0}}(\mathcal{O}_E, \mathcal{O}_E) 
\cong \bigwedge^p (\mathcal{O}_E(-1)^{m+1}).
\]
Hence
\[
\text{Hom}^p_{D(X)}(\mathcal{O}_E,\mathcal{O}_E) 
\cong \begin{cases} \mathbb{C} &\text{ if } p = 0,2m+1 \\
0 &\text{ otherwise}. 
\end{cases}
\]

(2) $\mathcal{O}_E$ in Example~\ref{flop}~(2) is {\it not} a $2m$-sherical 
object.
Indeed, since $N_{E/X_0} \cong \Omega^1_E$, we have
\[
\mathcal{E}xt^p_{\mathcal{O}_{X_0}}(\mathcal{O}_E, \mathcal{O}_E) 
\cong \Omega^p_E.
\]
Hence
\[
\text{Hom}^p_{D(X_0)}(\mathcal{O}_E,\mathcal{O}_E) 
\cong \begin{cases} \mathbb{C} &\text{ if } p = 0,2,\dots,2m \\
0 &\text{ otherwise}. 
\end{cases}
\]
\end{Expl}

There is some relationship between the flops and the twistings.

\begin{Expl}
If $m=1$ in Example~\ref{flop}~(1), then there are isomorphisms
\[
\Phi^{\mathcal{O}_Z}_{Y \to X} \circ \Phi^{\mathcal{O}_Z}_{X \to Y}
(\mathcal{O}_X(-k)) \cong 
T'_{\mathcal{O}_E(-1)}(\mathcal{O}_X(-k))
\]
for $k = 0,1$.

Indeed, we have 
\[
\begin{CD}
\mathcal{O}_X @>{\Phi^{\mathcal{O}_Z}_{X \to Y}}>> \mathcal{O}_Y
@>{\Phi^{\mathcal{O}_Z}_{Y \to X}}>> \mathcal{O}_X \\
\mathcal{O}_X(-1) @>{\Phi^{\mathcal{O}_Z}_{X \to Y}}>> \mathcal{O}_Y(1)
@>{\Phi^{\mathcal{O}_Z}_{Y \to X}}>> \mathcal{I}_E(-1)
\end{CD}
\]
where $\mathcal{I}_E$ is the ideal sheaf of $E$ in $X$.
On the other hand, 
\[
R\text{Hom}_X(\mathcal{O}_X,\mathcal{O}_E(-1))=0, \quad
R\text{Hom}_X(\mathcal{O}_X(-1),\mathcal{O}_E(-1))=\mathbb{C}
\]
hence 
\[
T'_{\mathcal{O}_E(-1)}(\mathcal{O}_X) = \mathcal{O}_X, \quad
T'_{\mathcal{O}_E(-1)}(\mathcal{O}_X(-1)) = \mathcal{I}_E(-1).
\]
\end{Expl}

\begin{Expl}
In Example~\ref{flop}~(2), there are isomorphisms
\[
\Phi^{\mathcal{O}_{Z_0}}_{X_0 \to Y_0} \circ
(\Phi^{\mathcal{O}_{\tilde Z_0}}_{X_0 \to Y_0})^{-1}(\mathcal{O}_{Y_0}(k))
\cong T_{\mathcal{O}_{F}(-1)}(\mathcal{O}_{Y_0}(k))
\]
for $k=0,1,\dots,m$.

Indeed, since 
\[
\begin{split}
&\text{Hom}_{D(Y_0)}^p(\mathcal{O}_{F}(-1), \mathcal{O}_{Y_0}(k))
\cong \text{Hom}_{D(Y_0)}(\mathcal{O}_{Y_0}(k), \mathcal{O}_{F}(-1)[2m-p])^* \\
&\cong \begin{cases} 0 &\text{ if } k=0,\dots,m-1 \\
0 &\text{ if } k=m \text{ and } p \ne m \\
\mathbb{C} &\text{ if } k=m \text{ and } p=m.
\end{cases}
\end{split}
\]
we have 
\[
T_{\mathcal{O}_{F}(-1)}(\mathcal{O}_{Y_0}(k)) \cong \mathcal{O}_{Y_0}(k)
\]
for $k=0,1,\dots,m-1$, and
\[
\mathcal{O}_{F}(-1)[-m] \to \mathcal{O}_{Y_0}(m) \to 
T_{\mathcal{O}_{F}(-1)}(\mathcal{O}_{Y_0}(m)) \to \mathcal{O}_{F}(-1)[-m+1] 
\]
is a distinguished triangle, where the first arrow is non-trivial.

On the other hand, we have an exact sequence
\[
0 \to \mathcal{O}_X \to \mathcal{O}_X \to \mathcal{O}_{X_0} \to 0 
\]
where the first arrow is the multiplication by an equation of $W_0 \subset W$.
Hence 
\[
\Phi^{\mathcal{O}_{\tilde Z_0}}_{X_0 \to Y_0}(\mathcal{O}_{X_0}(-k))
\cong \mathcal{O}_{Y_0}(k)
\]
for $k=0,1,\dots,m$.
If $k=0,1,\dots,m-1$, then we also have 
\[
\begin{CD}
\mathcal{O}_{X_0}(-k) @>{f_0^*}>> \mathcal{O}_{Z_0}(-k,0) \cong 
\mathcal{O}_{Z_0}(0,k)(kG_0) @>{g_{0*}}>> \mathcal{O}_{Y_0}(k)
\end{CD}
\]
because $K_{Z_0/Y_0} = (m-1)G_0$.

For $k=m$, we have an exact sequence
\[
0 \to \mathcal{O}_{Z_0}(0,m)((m-1)G_0) \to \mathcal{O}_{Z_0}(0,m)(mG_0)
\to \omega_{G_0}(0,m) \to 0.
\]
Since $g_{0*}(\omega_{G_0}) \cong \omega_{F}[-m+1] 
\cong \mathcal{O}_{F}(-m-1)[-m+1]$, 
we obtain a distinguished triangle 
\[
\mathcal{O}_{F}(-1)[-m] \to \mathcal{O}_{Y_0}(m) \to 
g_{0*}f_0^*(\mathcal{O}_{X_0}(-m)) \to \mathcal{O}_{F}(-1)[-m+1].
\]

We claim that the first arrow is non-trivial as an element of 
\[
\text{Hom}_{D(Y_0)}(\mathcal{O}_{F}(-1)[-m], \mathcal{O}_{Y_0}(m)) \cong 
\mathbb{C}.
\]
Indeed, if not, then we would have 
\[
\text{Hom}_{D(Y_0)}(\mathcal{O}_{F}(-1)[-m], \mathcal{O}_{Y_0}(m))
\cong \text{Hom}_{D(Y_0)}^1(g_{0*}f_0^*(\mathcal{O}_{X_0}(-m)), 
\mathcal{O}_{Y_0}(m))
\]
but 
\[
\begin{split}
&\text{Hom}_{D(Y_0)}^1(g_{0*}f_0^*(\mathcal{O}_{X_0}(-m)), 
\mathcal{O}_{Y_0}(m)) \\
&\cong \text{Hom}_{D(Z_0)}^1(f_0^*(\mathcal{O}_{X_0}(-m)), 
g_0^!\mathcal{O}_{Y_0}(m)) \\
&\cong \text{Hom}_{D(Z_0)}^1(\mathcal{O}_{Z_0}(0,m)(mG_0), 
\mathcal{O}_{Z_0}(0,m)((m-1)G_0)) \\
&\cong H^1(Z_0, \mathcal{O}_{Z_0}(-G_0)) \cong 0
\end{split}
\]
a contradiction.
Therefore, we have the desired isomorphisms.
\end{Expl}

\begin{Prop}\label{Mu2}
In Example~\ref{flop}~(2), if $m \ge 2$, then the functor
\[
g_{0*}f_0^* = \Phi^{\mathcal{O}_{Z_0}}_{X_0 \to Y_0}: D(X_0) \to D(Y_0)
\]
is not an equivalence.
\end{Prop}

\begin{proof}
Let us write $\Phi = \Phi^{\mathcal{O}_{Z_0}}_{X_0 \to Y_0}$ and 
$a = \mathcal{O}_{X_0}(-m)$. 
We consider a spectral sequence
\[
E_2^{p,q} = \bigoplus_{i \in \mathbb{Z}} \text{Ext}^p(H^i(\Phi(a)),
H^{q+i}(\Phi(a))) \Rightarrow \text{Hom}_{D(Y_0)}^{p+q}(\Phi(a),\Phi(a))
\]
given by the last line of \cite{V}~4.6.10.
We have
\[
H^q(\Phi(a)) \cong \begin{cases} \mathcal{O}_{Y_0}(m) &\text{ if } q = 0 \\
\mathcal{O}_{F}(-1) &\text{ if } q = m-1 \\
0 &\text{ otherwise} 
\end{cases}
\]
and
\[
\begin{split}
&\text{Ext}^p(\mathcal{O}_{Y_0}(m), \mathcal{O}_{Y_0}(m)) 
\cong \text{Ext}^p(\mathcal{O}_{X_0}(-m), \mathcal{O}_{X_0}(-m)) 
= 0 \text{ for } p \ne 0 \\
&\text{Ext}^p(\mathcal{O}_{Y_0}(m), \mathcal{O}_{F}(-1)) 
\cong \begin{cases} \mathbb{C} &\text{ if } p = m \\
0 &\text{ otherwise} \end{cases} \\
&\text{Ext}^p(\mathcal{O}_{F}(-1), \mathcal{O}_{Y_0}(m)) 
\cong \begin{cases} \mathbb{C} &\text{ if } p = m \\
0 &\text{ otherwise} \end{cases} \\
&\text{Ext}^p(\mathcal{O}_{F}(-1), \mathcal{O}_{F}(-1)) 
\cong \begin{cases} \mathbb{C} &\text{ if } p = 0,2,\dots, 2m \\
0 &\text{ otherwise.} \end{cases}
\end{split}
\]
Then the terms $E_2^{p,0}$ for $p = 2, \dots, 2m-2$ survive, 
hence $\text{Hom}_{D(X_0)}^p(a,a)$ and 
$\text{Hom}_{D(Y_0)}^p(\Phi(a),\Phi(a))$ are not isomorphic for these $p$.
\end{proof}

\begin{Rem}
After this paper was written, Jan Wierzba informed us that 
Corollary~\ref{Mu1} and Proposition~\ref{Mu2} were already proved by 
Namikawa~\cite{N}, though the proofs are different.
Combining with a result in \cite{CMSB} or \cite{WW} (see also \cite{Ke}), 
we obtain the implication from (2) to (1) in 
Conjecture~\ref{D and K}
in the case of symplectic projective manifolds of dimension $4$.
\end{Rem}


\section{Flops of terminal $3$-folds}

We shall deal with singular verieties in this section.

The smoothness of the given varieties is an important assumption for the study 
of derived categories.
For example, any coherent sheaf on a smooth projective variety has a finite 
locally free resolution, hence the Serre functor exists.

We can compare our situation with the deformation theory of maps 
from curves to varieties.
The latter is not applicable to singular varieties because 
the smoothness assumption is essential for a good obstruction theory.
However it provides deep results such as the theory of rationally
connected varieties.

We can still deal with singular varieties as if they are smooth in
some cases:
   
(1) If $X$ is a variety with only quotient singularities, then we consider 
a smooth stack $\mathcal{X}$ above $X$ 
as a natural substitute (cf. \cite{Francia}).

(2) If $X$ has only hypersurface singularities, then we embed $X$ into 
a smooth variety by deformations (cf. \cite{C}).

(3) If $X$ is a normal crossing variety, then we replace $X$ by its 
smooth hypercovering (cf. \cite{L}).

We consider a mixture of (1) and (2) in this section.

\begin{Defn}
Let $X$ be a normal quasiprojective 
variety such that the canonical divisor $K_X$ is a 
$\mathbb{Q}$-Cartier divisor.
Each point $x \in X$ has an open neighborhood $U_x$ such that $m_xK_X$ is a 
principal Cartier divisor on $U_x$ for a minimum positive integer $m_x$.
The {\it canonical covering} $\pi_x: \tilde U_x \to U_x$ is a finite 
morphism of degree $m_x$ from a normal variety 
which is etale in codimension $1$ and such that $K_{\tilde U_x}$ is a Cartier 
divisor. 
The canonical coverings are etale locally uniquely determined,
thus we can define the {\it canonical covering stack} $\mathcal{X}$ as 
the stack above $X$ given by the collection of canonical coverings 
$\pi_x: \tilde U_x \to U_x$.
\end{Defn}

We denote by $D(\mathcal{X}) = D^b(\text{Coh}(\mathcal{X}))$ 
the derived category of 
bounded complexes of coherent orbifold sheaves on $\mathcal{X}$
(cf. \cite{Francia}).

The following was suggested by Burt Totaro.

\begin{Prop}
Let $X$ be a normal projective 
variety such that the canonical divisor $K_X$ is a 
$\mathbb{Q}$-Cartier divisor.
Then there exists an embedding 
\[
\phi: X \to \mathbb{P}(a_1, \dots, a_N)
\]
to a weighted projective space such that 
the stack structure on $X$ induced from the natural smooth stack structure 
of $\mathbb{P}(a_1, \dots, a_N)$ coincides with the one
defined by the canonical coverings.
\end{Prop}

\begin{proof}
Let $H$ be an ample Cartier divisor such that $K_X + H$ is still ample as a 
$\mathbb{Q}$-Cartier divisor.
The ring $R = \bigoplus_{m=0}^{\infty}H^0(X, m(K_X+H))$ 
is a finitely generated algebra over $\mathbb{C}$.
Let $x_1, \dots, x_N$ be a set of homogeneous generators of $R$ of degree
$a_1, \dots, a_N$.
Then we obtain an embedding of $X$ to a weighted projective space
\[
\phi: X \to \mathbb{P}(a_1, \dots, a_N).
\]
Since $K_X + H$ is ample, $\text{g.c.d.}(a_1, \dots, a_N) = 1$.

We claim that 
\[
\text{g.c.d.}(a_1, \dots, \check{a_i}, \dots, a_N) = 1
\]
for any $i=1,\dots,N$, i.e., the sequence of integers 
$(a_1, \dots, a_N)$ is {\it well-formed}.
Indeed, suppose that 
$(a_2, \dots, a_N) = c \ne 1$.
Let $m$ be a sufficiently large integer which is not divisible by $c$, and 
consider an exact sequence
\[
\begin{split}
0 \to &\mathcal{O}_X((m-a_1)(K_X+H)) \to \mathcal{O}_X(m(K_X+H)) \\
&\to \mathcal{F}_m \to 0
\end{split}
\]
given by the multiplication by $x_1$, where $\mathcal{F}_m$ 
is a sheaf on $X_1=\text{div}(x_1)$.
By assumption, we have $H^0(X, (m-a_1)(K_X+H)) \cong H^0(X, m(K_X+H))$, while
$H^0(X_1, \mathcal{F}_m) \ne 0$ and $H^1(X, (m-a_1)(K_X+H)) = 0$
for large $m$, a contradiction. 

Let us fix a point $p \in X$. Then there exists a homogeneous coodinate,
say $x_1$, such that $x_1(p) \ne 0$.
We have a commutative diagram
\[
\begin{CD}
U @>{\phi}>> U_{x_1} @>{\subset}>> \mathbb{P}(a_1, \dots, a_N) \\
@A{\pi_U}AA   @AA{\pi_1}A  @. \\
\tilde U @>{\tilde{\phi}}>> \tilde U_{x_1} @.
\end{CD}
\]
where $U$ is a small open neighborhood of $p$,
$U_{x_1}$ is the open subset of $\mathbb{P}(a_1, \dots, a_N)$ defined
by $x_1 \ne 0$, 
$\pi_U: \tilde U \to U$ is a canonical covering, and 
$\pi_1: \tilde U_{x_1} \to U_{x_1}$ is the natural covering from 
an affine space with coordinates 
\[
x_2x_1^{-a_2/a_1}, \dots, x_Nx_1^{-a_N/a_1}.
\]
Note that both $\pi_U$ and $\pi_1$ are etale in codimension $1$.

Since $x_1(p) \ne 0$, we may choose a branch of $x_1^{1/a_1}$ on
sufficiently small $U$.
Then $\phi$ can be lifted to a morphism  
$\tilde{\phi}: \tilde U \to \tilde U_{x_1}$ 
which we can check to be etale.
Therefore, the two stack structures coincide.
\end{proof}

\begin{Rem}
(1) By the proposition, any coherent orbifold sheaf on the canonical 
covering stack $\mathcal{X}$ 
has a surjection from a locally free orbifold sheaf on $\mathcal{X}$.
But the Serre functor for the category $D(\mathcal{X})$ does not exist 
in general.

(2) Totaro (\cite{T}) proved the following resolution theorem:
on a smooth orbifold whose coarse moduli space is a separated scheme,
any coherent orbifold sheaf has a finite resolution by 
locally free orbifold sheaves.
\end{Rem}

We still have a good spanning class for terminal $3$-folds:

\begin{Lem}\label{span}
Let $X$ be a normal projective variety of dimension $3$ 
with only terminal singularities, $m_x$ the index of $K_X$ at $x \in X$,
and $\mathcal{X}$ the canonical covering stack of $X$.  
Then the set $\{\mathcal{O}_x(iK_X) \vert x \in X, 0 \le i < m_x\}$ 
is a spanning class of $D(\mathcal{X})$.
\end{Lem}

\begin{proof}(cf. \cite{B1}~Example~2.2 and \cite{C}~Lemma~3.4)
Let $a$ be a non-zero object of $D(\mathcal{X})$.
Take a point $x_0$ in the support of $a$, and let
$q_0$ be the maximal value of $q$ such that $H^q(a)_{x_0} \ne 0$.
Then there exists an integer
$i_0$ such that $\text{Hom}(H^{q_0}(a), \mathcal{O}_{x_0}(i_0K_X)) \ne 0$.
Then $\text{Hom}^{-q_0}_{D(X)}(a, \mathcal{O}_{x_0}(i_0K_X)) \ne 0$.

If the support of $a$ is not contained in the singular locus of $X$, then 
we take the above point $x_0$ from the smooth locus of $X$.
By the Serre duality, we have $\text{Hom}^{n+q_0}(\mathcal{O}_{x_0}, a) \ne 0$,
where $n = \dim X$.
Otherwise, let $q_1$ be the minimal value of $q$ such that $H^q(a)_{x_0} 
\ne 0$. 
Since $X$ has only isolated singularities, there exists an integer
$i_1$ such that $\text{Hom}(\mathcal{O}_{x_0}(i_1K_X), H^{q_1}(a)) \ne 0$.
Hence $\text{Hom}^{q_1}_{D(X)}(\mathcal{O}_{x_0}(i_1K_X), a) \ne 0$.
\end{proof}

\begin{Thm}\label{flop to D}
Let $X$ and $Y$ be normal quasiprojective varieties of dimension $3$ 
with only $\mathbb{Q}$-factorial terminal singularities,  
\[
\begin{CD}
X @>{\phi}>> W @<{\psi}<< Y
\end{CD}
\]
a flop, and $\mathcal{X}$ and $\mathcal{Y}$ the canonical covering stacks 
above $X$ and $Y$, respectively.
Then the bounded derived categories of coherent orbifold sheaves 
$D(\mathcal{X})$ and $D(\mathcal{Y})$ are equivalent
as triangulated categories.
\end{Thm}

\begin{proof}
The assertion is already proved in the case where $K_X$ is a Cartier divisor
by Bridgeland~\cite{B2} and Chen~\cite{C} (see also \cite{VB}).
Indeed, it is proved that the structure sheaf $\mathcal{O}_Z$ 
of the fiber product $Z = X \times_W Y$ is quasi-isomorphic to 
a finite complex of sheaves on $X \times Y$ flat over $X$ so that the 
integral functor $\Phi^{\mathcal{O}_Z}_{X \to Y}: D(X) \to D(Y)$
is defined and is an equivalence (\cite{C}~Lemma~2.1 and 
Proposition~4.2).

We shall give a new simpler proof, which is based on \cite{VB}~\S 4.1, that 
$\Phi = \Phi^{\mathcal{O}_Z}_{X \to Y}: D(X) \to D(Y)$ is an equivalence in
the case where $K_X$ is a Cartier divisor.
We may assume that $W$ is a hypersurface singularity of multiplicity $2$.
Thus $W$ has an involution $\sigma$ such that $W/\langle \sigma \rangle$ is 
smooth.
We may take $Y = X$ and $\psi = \sigma \circ \phi$.  

First we prove that $\Phi(\mathcal{O}_X) \cong 
\mathcal{O}_Y$. 
Indeed, for any closed point $y \in Y$, the scheme theoretic fiber 
$g^{-1}(y)$ is isomorphic to the fiber $f^{-1}(\psi(y))$.
We have $H^0(\mathcal{O}_{f^{-1}(\psi(y))}) \cong \mathbb{C}$, 
hence the natural homomorphism 
$\mathcal{O}_Y \to R^0g_*\mathcal{O}_Z$ is an isomorphism.

Any subscheme of $Z$ which is mapped by $g$
to an infinitesimal subscheme $\bar y$ of $Y$ supported at $y$  
is isomorphic to a subscheme of the product $\bar C \times \bar y$ 
for a subscheme $\bar C$ of $X$ which is mapped by $\phi$
to an infinitesimal subscheme $\bar w$ of $W$ supported at $\phi(y)$.
Since $R^1\phi_*\mathcal{O}_X = 0$, it follows that $R^1g_*\mathcal{O}_Z = 0$.

Let $C_j$ ($j=1, \dots, t$) be the exceptional curves of $\phi$, and
$L_i$ ($i=1, \dots, t$) invertible sheaves on $X$ such that
$(L_i \cdot C_j) = \delta_{ij}$.
Then $L_i$ are generated by global sections for all $i$.
We note that $R^1\phi_*L_i^*$ may not necessarily vanish.

According to \cite{VB}~\S 4.1, we construct 
locally free sheaves $M_i$ and $N_i$ on $X$ 
by the following exact sequences
\[
\begin{split}
&0 \to \mathcal{O}_X^{r_i} \to M_i \to L_i \to 0 \\
&0 \to N_i \to \mathcal{O}_X^{s_i} \to L_i \to 0 
\end{split}
\]
for some integers $r_i, s_i$ such that we have
the vanishing higher direct image sheaves 
$R^1\phi_*M_i^* = 0$ and $R^1\phi_*N_i = 0$.
By \cite{VB}~Proposition~4.1.2, if we take $r_i$ and $s_i$ to be the minimal
possible integers under the vanishing conditions, then 
we have 
\[
\phi_*N_i \cong \sigma_*\phi_*M_i
\]
where we note that $\sigma_*L_i \cong L_i^*$.
By construction, $M_i$ and $N_i^*$ are generated by global sections.

It follows that $R^1g_*f^*N_i = 0$ from $R^1\phi_*N_i = 0$ as before.  
We consider an exact sequence
\[
0 \to g_*f^*N_i \to \mathcal{O}_Y^{s_i} \to g_*f^*L_i \to 0.
\]
Since there is a non-natural injection $g_*f^*L_i \to g_*\mathcal{O}_Z$, 
the sheaf $g_*f^*L_i$ is torsion free.
Hence $g_*f^*N_i$ is a reflexive sheaf.
Since $\psi_*g_*f^*N_i \cong \phi_*N_i \cong \psi_*M_i$,
we conclude that $\Phi(N_i) \cong M_i$.

The set of sheaves $\Omega = \{\mathcal{O}_X, N_1, \dots, N_t\}$ 
is a spanning class of $D(X)$.
$\omega$ is locally free, $\omega^*$ is generated by global sections and 
$R^1\phi_*\omega = 0$ for any $\omega \in \Omega$. 
Hence 
\[
\text{Hom}_{D(X)}^p(\omega_1, \omega_2) = 0
\]
for $p > 0$ and $\omega_1, \omega_2 \in \Omega$.
Similarly we have
\[
\text{Hom}_{D(Y)}^p(\Phi(\omega_1), \Phi(\omega_2)) = 0.
\]
Since $X$ and $Y$ are isomorphic in codimension $1$, we have 
\[
\text{Hom}_X(\omega_1, \omega_2)
\cong \text{Hom}_Y(\Phi(\omega_1), \Phi(\omega_2)).
\]
Therefore, we have proved that $\Phi$ is an equivalence in the case where 
$K_X$ is a Cartier divisor.

Now we consider the general case.
Let $\mathcal{W}$ be the canonical covering stack of $W$.
Let $w \in W$ be a point, $W_w$ its small neighborhood 
on which $m_wK_W$ is a principal Cartier divisor, and 
$\pi_w: \tilde W_w \to W_w$ a canonical covering.
Then $m_wK_X$ and $m_wK_Y$ are also principal Cartier divisors on
$X_w = \phi^{-1}(W_w)$ and $Y_w = \psi^{-1}(W_w)$, respectively, and
we have corresponding canonical coverings
$\pi_X: \tilde X_w \to X_w$ and $\pi_Y: \tilde Y_w \to Y_w$.
Thus there are morphisms of stacks $\phi: \mathcal{X} \to \mathcal{W}$
and $\psi: \mathcal{Y} \to \mathcal{W}$.
Let 
\[
\mathcal{Z} = \mathcal{X} \times_{\mathcal{W}} \mathcal{Y}
\] 
be the fiber product as a stack.
Then it is a stack above $Z = X \times_W Y$ where 
local coverings are given by
\[
\tilde Z_w = \tilde X_w \times_{\tilde W_w} \tilde Y_w 
\to Z_w = X_w \times_{W_w} Y_w.
\]
Let $\mathfrak{f}: \mathcal{Z} \to \mathcal{X}$
and $\mathfrak{g}: \mathcal{Z} \to \mathcal{Y}$ be the induced morphisms.

We claim that the functor 
\[
\mathfrak{g}_*\mathfrak{f}^*: 
D(\mathcal{X}) \to D(\mathcal{Y})
\]
is defined and is an equivalence.
Indeed, over an open subset $W_w$, we know already that the integral functor
\[
\Phi^{\mathcal{O}_{\tilde Z_w}}_{\tilde X_w \to \tilde Y_w}: 
D(\tilde X_w) \to D(\tilde Y_w)
\]
is an equivalence.
Let $\mathcal{X}_w = \mathcal{X} \vert_{X_w} = [\tilde X_w/G]$,
$\mathcal{Y}_w = \mathcal{Y} \vert_{Y_w} = [\tilde Y_w/G]$, 
$\mathcal{Z}_w = \mathcal{Z} \vert_{Z_w} = [\tilde Z_w/G]$, 
$\mathfrak{f}_w = \mathfrak{f} \vert_{\mathcal{Z}_w}$ and 
$\mathfrak{g}_w = \mathfrak{g} \vert_{\mathcal{Z}_w}$.
The Galois group $G = \mathbb{Z}/m_w$ acts equivariantly so that we have
$D(\tilde X_w)^G \cong D(\mathcal{X}_w)$ and 
$D(\tilde Y_w)^G \cong D(\mathcal{X}_w)$ (cf. \cite{BM2}).
Hence we have a well-defined equivalence
\[
\mathfrak{g}_{w*}\mathfrak{f}_w^*: 
D(\mathcal{X}_w) \to D(\mathcal{Y}_w).
\]
By Lemma~\ref{span}, we conclude the proof.
\end{proof}

\begin{Rem}
We note that the equivalence 
$\Phi=\mathfrak{g}_*\mathfrak{f}^*: D(\mathcal{X}) \to D(\mathcal{Y})$ does not
induce an equivalence $D(X) \to D(Y)$ of usual derived categories for
singular varieties.
Indeed, we can construct a similar example as in \cite{Francia}~Example~5.1.
There is a skyscraper sheaf $a \in D(\mathcal{X})$ supported over
a non-Gorenstein singular point of $X$ such that $\pi_{X*}(a) = 0$ in $D(X)$,
but its image $\Phi(a) \in D(\mathcal{Y})$ has a $1$-dimensional support so
that $\pi_{Y*}(\Phi(a)) \ne 0$ in $D(Y)$.
\end{Rem}


Department of Mathematical Sciences, University of Tokyo, 

Komaba, Meguro, Tokyo, 153-8914, Japan 

kawamata@ms.u-tokyo.ac.jp

\end{document}